\documentclass[11pt]{article}

\usepackage{mathrsfs}
\usepackage{amsfonts}
\usepackage[leqno]{amsmath}
\usepackage{graphicx}
\usepackage{latexsym}
\usepackage{float}

\usepackage{amsmath,amsfonts,amssymb,amsthm,mathrsfs,euscript,
	color}
\usepackage{enumerate}

\def\sqr#1#2{{\vcenter{\vbox{\hrule height.#2pt
				\hbox{\vrule width.#2pt height#1pt \kern#1pt \vrule width.#2pt}
				\hrule height.#2pt}}}}
\def\5n{\negthinspace \negthinspace \negthinspace \negthinspace \negthinspace }
\def\4n{\negthinspace \negthinspace \negthinspace \negthinspace }
\def\3n{\negthinspace \negthinspace \negthinspace }
\def\2n{\negthinspace \negthinspace }
\def\1n{\negthinspace }

\def\EE{\mathsf E\:\!}
\def\PP{\mathsf P}

\newcommand{\parens}[1]{\left(#1\right)}
\newcommand{\bracks}[1]{\left[#1\right]}

\def\dbR{\mathbb{R}}

\def\sC{\mathscr{C}}

\def\cF{{\cal F}}

\def\cT{{\cal T}}

\def\cX{{\cal X}}

\def\BD{{\bf D}}

\def\BS{{\bf S}}

\def\ds{\displaystyle}

\def\ns{\noalign{\smallskip}}

\def\q{\quad}
\def\qq{\qquad}

\def\({\Big (}
\def\){\Big )}
\def\[{\Big[}
\def\]{\Big]}

\def\e{\varepsilon}

\def\k{\kappa}

\def\t{\tau}

\def\th{\theta}

\def\wt{\widetilde}

\def\cd{\cdot}

\def\square#1{\vbox{\hrule\hbox{\vrule height#1%
			\kern#1\vrule}\hrule}}
\def\rectangle#1#2{\vbox{\hrule\hbox{\vrule height#1%
			\kern#2\vrule}\hrule}}

\font\tenbb=msbm10 \font\sevenbb=msbm7 \font\fivebb=msbm5

\newfam\bbfam
\scriptscriptfont\bbfam=\fivebb \textfont\bbfam=\tenbb
\scriptfont\bbfam=\sevenbb

\newtheorem{theorem}{\indent Theorem}[section]

\newtheorem{proposition}[theorem]{\indent Proposition}

\newtheorem{lemma}[theorem]{\indent Lemma}
\newtheorem{remark}[theorem]{\indent Remark}

\newtheorem{assumption}[theorem]{\indent Assumption}

\makeatletter

\@addtoreset{equation}{section}
\makeatother

\def\bea{\begin{equation*}}
\def\eea{\end{equation*}}

\def\bel{\begin{equation}\label}
\def\eel{\end{equation}}

\def\ba{\begin{array}}
	\def\ea{\end{array}}
\newcommand{\ad}{&\!\!\!\displaystyle}

\def\({\Big (}
\def\){\Big )}
\def\[{\Big[}
\def\]{\Big]}
\def\q{\quad}
\def\qq{\qquad}

\def\e{\varepsilon}

\def\ds{\displaystyle}

\def\ns{\noalign{\smallskip}}

\def \hat {\widehat}
\def\LL{I\!\!L}

\def\para#1{\vskip .25\baselineskip\noindent{\bf #1}}
\allowdisplaybreaks[4]
\begin{document}

	\title{ The Minimax Wiener Sequential Testing Problem}
	
	\author{Philip A. Ernst\footnote{     Department of Mathematics, Imperial College London, London SW7 2AZ, United Kingdom. Email: {\tt p.ernst@imperial.ac.uk}.} ~~and ~~Hongwei Mei\footnote{ Department of Mathematics and Statistics, Texas Tech University, Lubbock, TX79409.   Email: {\tt hongwei.mei@ttu.edu.} } }
	\maketitle	
	\begin{abstract}
		
			Consider the sample path of a one-dimensional diffusion for which the diffusion coefficient is given and where the drift may take on one of two values: $
\mu_0$ or $\mu_1$. Suppose that the signal-to-noise ratio (defined as the difference between the two 
possible   drifts 
   divided by the diffusion coefficient) is non-constant.  Given an initial state for the observed process, we consider a minimax formulation of the Wiener sequential testing problem for detecting the correct drift coefficient as soon as possible and with minimal probabilities of incorrect terminal decisions. We solve the problem in the Bayesian formulation, under any prior probabilities of the process having drift $\mu_0$ or $\mu_1$, when the passage of time is
penalized linearly. In the case where the signal-to-noise ratio is assumed constant,  we obtain an explicit formula for the least favorable distribution.

	\end{abstract}

	\para{Keywords}: Wiener sequential testing, minimax optimization, least favorable distribution, optimal stopping, stochastic differential equations.
	
	\para{MSC 2010 Codes}: Primary: 60G40, 93C30. Secondary: 60H30, 91B70.

	\maketitle

	\section{Introduction}\label{sec1}

	Consider the real time observation of a one-dimensional diffusion process $X$ whose dynamics are given by the following stochastic differential equation (SDE)
	\bel{SDE0} dX_t=\big(\mu_0(X_t)+\th(\mu_1(X_t)-\mu_0(X_t)\big )dt+\sigma(X_t)dB_t,\eel 
 where  $B$ is a one dimensional standard Brownian motion and $\mu_0$, $\mu_1$  and $\sigma$ are appropriately defined functions assumed as given (see Assumption \ref{a1}). The unobservable random variable $\th$ denotes the true drift coefficient of $X$; i.e. if $\th=0$ ($\th=1$ correspondingly), the drift coefficient of $X$ is $\mu_0$ ($\mu_1$ correspondingly). Given that the process $X$ is observed in real time, the problem is to detect the correct drift coefficient as soon as possible and with minimal probabilities of incorrect terminal decisions.  \\ 
\indent The above sequential testing problem of two simple hypotheses for the drift of a one-dimensional diffusion is known as the \textit{Wiener sequential testing problem} (\cite{Gape2004,PS}). The Wiener sequential testing problem admits two distinct formulations. In the Bayesian framework, the random variable $\th$ takes on the value 0 with prior probability $1-\pi$ and the value 1 with prior probability $\pi$ (see \cite{Wald}). The variational (frequentist) paradigm makes no a priori probabilistic assumptions about $\theta$. In the sequel, we shall exclusively consider the Bayesian formulation.\\
\indent There is a well-established body of literature on the Wiener sequential testing problem. The books by Peskir and Shiryaev (\cite[Chapter IV]{PS}) and by Shiryaev (\cite[Chapter IV]{Sh-2}) solve the one-dimensional  Wiener sequential testing problem for an infinite horizon and Gapeev and Peskir (\cite{Gape2004}) solve the problem in the case of a finite horizon. The Wiener sequential testing problem was also recently considered in two and three dimensions by Ernst, Peskir, and Zhou (\cite{ErPZ2020}). A critical assumption in all the aforementioned references is that the {\it signal-to-noise ratio} (defined as the
difference between the two possible drifts divided by the diffusion coefficient) is constant. In the case of constant signal-to-noise ratio (SNR), the Wiener sequential testing problem is linear because the standard likelihood ratio process satisfies a linear stochastic differential equation
(\cite[Chapter IV]{PS}). If the signal-to-noise ratio is no longer assumed to be constant, the Wiener sequential testing problem becomes nonlinear.\\
\indent 
The nonlinear Wiener sequential testing problem was first considered by Gapeev and Shiryaev (\cite{Gape2011}), who solved the problem under the assumption of the existence of a regular solution to a free boundary problem. The more general setting of a non-constant signal-to-noise ratio has also been considered in \cite{JP2018}, where $X$ is assumed to be a Bessel process of dimension 2 or 3, as well as in \cite{JP2021}, where it is assumed that the drift can be distributed according to a known probability law. In both \cite{JP2018} and \cite{JP2021}, the optimal stopping boundaries
are characterized as the unique solution to a coupled system of nonlinear Volterra integral
equations and the solutions to the corresponding optimal stopping problems are functions of $\pi$ (cf. \cite[Equation (9.11)]{JP2018} or \cite[Equation (8.20)]{JP2021}).

\indent  This paper is in part motivated by the following question: is it possible (within the Bayesian setting) to reformulate the statement of the nonlinear Wiener sequential testing problem so that the solution to the optimal stopping problem does not depend on knowledge of the a priori distribution $\pi$? The answer is in the affirmative when one considers a minimax version of the Wiener sequential testing problem. In the minimax optimization framework, the objective is to minimize the cost functional under the worst case scenario of all possible prior distributions. The prior distribution of $\pi$ which yields the worst scenario among all possible a priori distributions is called the \textit{least favorable distribution} (see, for example, \cite{Lehmann} and references therein). Providing a characterization of this least favorable prior distribution solves the minimax optimization problem.  To the best of our knowledge, the present paper is the first to consider a minimax optimization framework for either the linear Wiener sequential testing problem or for the nonlinear Wiener sequential testing problem.  \\
\indent On the way to solving the minimax Wiener sequential testing problem in Section \ref{sec:mst}, we shall find it necessary to return to the work of Gapeev and Shiryaev (\cite{Gape2011}), who solved the nonlinear Wiener sequential testing problem under the assumption of the existence of unique solution to a free-boundary problem. Section \ref{sec:osp} of the present paper provides a solution to the nonlinear Wiener sequential testing problem under much weaker assumptions (see Assumption \ref{time-changeok} and Assumption \ref{time-changeok-2}). Our assumptions shall only concern the values of the two possible drift coefficients ($\mu_0$ and $\mu_1$) and the diffusion coefficient $\sigma$. This is the first contribution we make in Section \ref{sec:osp}. Our second contribution in Section \ref{sec:osp} is the proof of probabilistic regularity (cf. \cite[Definition 2.9, p.245]{KA1998}) of the corresponding optimal stopping boundaries. These contributions are encapsulated by Theorem \ref{ospt}. The key mathematical difficulties in this part of the paper arise from being unable to assume the existence of a unique solution to the free-boundary problem.\\
\indent The remainder of the paper is organized as follows. In Section \ref{sec:for}, we provide the formulation of two optimal stopping problems under consideration: (i) the nonlinear Wiener sequential testing problem (ii) the minimax nonlinear Wiener sequential testing problem. Section \ref{sec:osp} is devoted to providing a solution to the nonlinear Wiener sequential testing problem without invoking the strong assumption of the existence of a unique solution to the free-boundary problem (\cite{Gape2011}). In Section \ref{sec:mst}, we solve the minimax nonlinear Wiener sequential testing problem 
by finding a characterization of the least favorable distribution. Section \ref{sec:SNRC} is concerned with the minimax Wiener sequential testing problem in the special case where SNR is assumed constant. It is in this particular case that we are able to obtain an explicit formula for the least favorable distribution.

\section{Formulation of the optimal stopping problems}\label{sec:for}
\indent We begin with some necessary notation. For the sake of consistency, we shall follow the same notation used in Gapeev and Shiryaev (\cite{Gape2011}). \smallskip

1.
Given a probability space $(\Omega, \cF, \PP^\psi)$, let the measure $\PP^\psi$ be defined, for any $\psi\in [0,\infty]$, as 
\begin{align}\label{defPpi}
\PP^\psi(\cdot)=\frac1{1+\psi}\PP(\cdot|\th=0)+\frac\psi{1+\psi}\PP(\cdot|\th=1),
\end{align}
where  $\th$ is a random variable taking values in $\{0,1\}$ with $\PP^\psi(\th=1)=\psi/(1+\psi)$ and $\PP^\psi(\th=0)=1/(1+\psi)$. Let the measure $\PP^0$ be defined as $\PP^0(\cdot)=\PP(\cdot|\th=0)$ and let the measure $\PP^\infty$ be defined as $\PP^\infty(\cdot)=\PP(\cdot|\th=1)$.\\ 
\indent Suppose that we observe a process $X$ whose dynamics are given by the solution to the SDE in \eqref{SDE0}. With the exception of Section \ref{sec:SNRC}, we shall assume throughout the present paper that the signal-to-noise ratio, defined as
\bel{SNR}\rho(x):=\frac{(\mu_1(x)-\mu_0(x))}{\sigma(x)},\eel
is non-constant.  We shall also assume that Assumption \ref{a1} below holds throughout the sequel.

\begin{assumption}\label{a1} Let $\cX$ be a domain in $\dbR$.

{\rm (1).} The functions $\mu_i(\cdot)$ and $\sigma(\cdot)$ are continuously differentiable in $\cX$. Further, for any $x\in \cX$, $\sigma(x)>0$. \smallskip

	{\rm (2).} 
	For either $\th=0$ or $\th=1$, the SDE \eqref{SDE0} admits a unique solution in $\cX$. Moreover, $\int_0^t\rho^2(X_s)ds$ strictly increases  to $\infty$ almost surely as $t\rightarrow\infty$.
\end{assumption}
\noindent Part (1) of Assumption \ref{a1} offers a slight generalization of the work of Gapeev and Shirayev (\cite{Gape2011}), where $\cX$ is assumed to be $\cX=(0,\infty)$.
Part (2) of Assumption \ref{a1} will be invoked to guarantee that the likelihood ratio will either tend to infinity or tend to 0 when the waiting time tends toward infinity.
\smallskip

2. 
Being based upon the continuous observation of $X$, the problem is to sequentially test the two
hypotheses $H_0$ and $H_1$ with minimal loss where
\begin{align*}H_0:\th=0\qq\text{ and }\qq H_1:\th=1.
\end{align*}
We achieve this task by considering a sequential decision rule
$(\tau, d)$, where $\tau$ is a stopping time of the observation process $X$ (i.e. a stopping time with respect to the natural filtration $\cF^X_t=\sigma(X_s:0\leq s\leq t)$ for $t\geq 0$) and $d$ is an $\cF_{\tau}$-measurable function taking values in $\{0,1\}$. After we stop observing the process at time $\tau$, the decision function $d$ delineates which hypothesis should be accepted: if $d=0$, we accept $H_0$ and if $d=1$, we accept $H_1.$ The set of all admissible decisions can be written as $${\bf D}=\Big\{(\tau,d): \tau\text{ is a stopping time of $X$ and $d$ is a map valued in $\{0,1\}$}\Big\}.$$ 
\indent Our objective now is to find the optimal $(\tau^*,d^*)$ which minimizes the risk function 
\bel{per-1}J(x,\psi;\tau,d)
:=\PP^{\psi}(\th=0,d=1)+\PP^{\psi}(\th=1,d=0)+c\EE^\psi\tau,\eel
where $ X_0=x$ and  $c$ is a given constant. In the above equation, the combination $\PP^{\psi}(\th=0,d=1)+\PP^{\psi}(\th=1,d=0)$ 
is the {\it probability of wrong detection }and the term $\EE^\psi\tau$ is the {\it expected waiting time}.
\smallskip

3. As in \cite{Gape2011, PS}, we solve the equivalent optimal stopping problem which arises after performing an appropriate change of measure to the observed process $X$. 
A key observation in both \cite{Gape2011} and \cite{PS} is that
the optimal decision at time $\tau$ satisfies
\bel{decisionprocedure}\bar d_{\tau}=\left\{\ba{ll}0,\qq\text{ if }\psi L_{\tau}<1\\1,\qq\text{ if }\psi L_{\tau}> 1,\\
\text{either 0 or 1}, ~~\text{ if }\psi L_{\tau}=1,\ea\right.\eel
where the {\it likelihood ratio} $L_t$ is defined as  \bel{liklihood}L_t=\frac{d\PP^\infty}{d\PP^0}|_{\cF^X_{t}}.\eel 
The form of the optimal decision in \eqref{decisionprocedure} results from minimizing the probability of false detection given by \eqref{per-1}. The class of admissible decisions $\BD$ reduces to 
$${\bf S}=\{\tau: \tau\text{ is a stopping time of }X \}.$$
The set {\bf S} is thus to be interpreted as the subset of {\bf D} for which \eqref{decisionprocedure} holds. The expression in \eqref{decisionprocedure} enables us to rewrite the risk function in \eqref{per-1} as
\bel{barJ}\bar J(x,\psi;\tau)=J(x,\psi;\tau,\bar d_\tau).\eel
\indent The nonlinear Wiener sequential testing problem may now be formulated as the following optimal stopping problem (OSP). \smallskip

\noindent {\bf OSP 1}: Find a $\tau^*\in \BS$ such that 
\bel{osp}\bar J(x,\psi;\tau^*)=\inf_{\tau\in{\bf D}}\bar J(x,\psi;\tau).\eel
As noted in Section \ref{sec1}, the optimal stopping problem in \eqref{osp} was first solved by Gapeev and Shiryaev (\cite{Gape2011}). It is proven   (see Lemma 3.1 therein) that the optimal stopping time for this optimal stopping problem satisfies 
\bel{optimalstoppingtime} \tau^*=\inf\Big\{t\geq 0:  \psi L_t\notin (l^*_0(X_t),l^*_1(X_t))\Big\},\eel
for some appropriate functions $l^*_0(\cdot), l^*_1(\cdot)$ independent of $\psi$. We refer the reader to \cite{Gape2011} for further details. 
\smallskip

4. With the above preparation in hand, we turn to the formulation of the minimax Wiener sequential testing problem. The objective of the minimax formulation is to minimize the performance functional $\bar J(x,\psi;\cdot)$ in the \textit{worst case scenario} of all prior distributions $\psi$.  This leads to the following optimal stopping problem. \smallskip

\noindent {\bf OSP 2}: Find a $\tau^*\in{\bf S}$ such that  \bel{minimaxV1}\sup_{\psi\geq 0}\bar J(x,\psi;\tau^*)=\inf_{\tau\in{\bf S}}\sup_{\psi\geq0}\bar J(x,\psi;\tau).\eel
In the above equality, the ``worst case scenario" for $\psi$ corresponds to taking the supremum (over all $\psi \geq 0$) of  $\bar{J}(x,\psi;\tau)$.
\smallskip

5. In order to solve the optimal stopping problem in \eqref{minimaxV1}
by the saddle point property, we will need to find an optimal couple $(\varphi_0,\tau_0)$  satisfying,  for any $\psi\geq 0$ and $\tau\in{\bf S}$,
\begin{equation} \label{inequality} 
\bar J(x_0,\psi;\tau_0)\leq \bar J(x_0,\varphi_0;\tau_0)\leq \bar J(x_0,\varphi_0;\tau).
\end{equation}
The first inequality in \eqref{inequality} states that $\varphi_0$ is the least favorable  priori distribution given the stopping time $\tau_0$. The second inequality in \eqref{inequality} asserts that $\tau_0$ is the solution for the optimal stopping problem with initial value $\psi_0$.  By \eqref{optimalstoppingtime},  it follows that 
$$\tau_0=\inf\{t\geq 0:L_t\notin (l_0^*(X_t^{x_0})/\varphi_0,(l_1^*(X_t^{x_0})/\varphi_0\}=:\tau^*(\varphi_0),$$
where $l_0^*(\cdot)$ and $l_1^*(\cdot)$ are determined by the optimal stopping boundary in \eqref{osp}.\\
\indent This paper's key result for the minimax Wiener sequential testing problem, given by Theorem \ref{thm41}, is that there exists a $\varphi_0$ which is the least favorable distribution for the stopping time $\tau^*(\varphi_0)$. However, before embarking on this problem, we must first return to the nonlinear Wiener sequential testing problem in \eqref{osp}.

\section{The optimal stopping problem in \eqref{osp}}\label{sec:osp}
This section is concerned with the nonlinear Wiener sequential testing problem in \eqref{osp}. This problem was first solved by Gapeev and Shiryaev (\cite{Gape2011}) under the strong assumption of the existence of a unique solution to the free-boundary problem. Our interest in Theorem \ref{ospt} is to solve the same problem under much weaker assumptions. As we shall see, this in turn enables us to solve the minimax Wiener sequential testing problem in \eqref{minimaxV1}.\\
\indent Our solution to the optimal stopping problem in \eqref{osp} shall rely on one of two significantly weaker assumptions. These assumptions are given by Assumption \ref{time-changeok} and Assumption \ref{time-changeok-2} below. Note that each assumption only depends on the values of the two possible drift coefficients ($\mu_0$ and $\mu_1$) and the diffusion coefficient $\sigma$.

\begin{assumption}\label{time-changeok} Suppose  $\rho^2(\cdot)$  is decreasing on $\cX$  and, for all $x\in\cX$,
 \bel{DefK}K(x):= \frac{\mu_0(x)}{\mu_1(x)-\mu_0(x)}-\frac12\parens{\frac{\sigma^2(x)}{\mu_1(x)-\mu_0(x)}}'>-\frac12.\eel
\end{assumption} 

\begin{assumption}\label{time-changeok-2} Suppose $\rho^2(\cdot)$  is increasing on $\cX$ and, for all $x\in\cX$,
\bel{DefK2}K(x)<-\frac12.\eel
\end{assumption} 
\noindent Assumptions \ref{DefK} and \ref{DefK2} shall also be employed to verify that the optimal stopping boundary for the optimal stopping problem in \eqref{osp} is \textit{probabilistically regular} (cf. \cite[p.245]{KA1998}).  \\
 \indent We now proceed to express an equivalent formulation of the optimal stopping problem in \eqref{osp}.  
We consider the {\it likelihood ratio process}  $L$ defined by
\bel{liklihood2}L_t=\frac{d\PP^\infty}{d\PP^0}|_{\cF^X_{t}}.\eel
Invoking Girsanov's theorem, we calculate
\bel{Lt}L_t=\frac{\exp\parens{\ds \int_0^t\langle\sigma^{-1}(X_s)\mu_1(X_s),dX_s\rangle-\frac12\int_0^t\langle \sigma^{-1}(X_s)\mu_1(X_s),\mu_1(X_s)\rangle ds}}{\exp\parens{\ds\int_0^t\langle\sigma^{-1}(X_s)\mu_0(X_s),dX_s\rangle-\frac12\int_0^t\langle \sigma^{-1}(X_s)\mu_0(X_s),\mu_0(X_s)\rangle ds}}.\eel
It may be easily checked that $L$ satisfies
\bel{Lt2}\frac{dL_t}{L_t}=\rho(X_t) \parens{dX_t-\mu_0(X_t)dt},\text{ with }L_0=1.\eel 
We proceed to apply a change of measure on the performance functional $J$ given in \eqref{barJ}. Applying Girsanov's theorem yields
$$\frac{d\PP^0}{d\PP^\psi}\Big|_{\cF^X_t}=\frac{1+\psi L_t} {1+\psi}.$$
We pause to note that
$$ \ba{ll}\ad \bar J(x,\psi;\tau)=\EE^\psi_{x}\bracks{\frac{1}{1+\psi L_\tau}\wedge \frac{\psi L_\tau }{1+\psi L_\tau}+c\int_0^\tau 1dt}\\
\ns\ad=\frac{1}{1+\psi}\EE^0_{x}\bracks{1\wedge (\psi L_t)+c\int_0^\tau ( 1+\psi L_t)dt}.\ea$$
Let us define $\varPhi_t=\psi L_t$. Note that 
\bel{SDE}\left\{\ba{ll}\ad dX_t=\mu_0(X_t)dt+\sigma(X_t)d\wt B_t\\
\ns\ad  d\varPhi_t=\rho(X_t)\varPhi_t d\wt B_t,\ea\right.\eel
where $$d\wt B_t=\sigma^{-1}(X_t)\parens{dX_t-\mu_0(X_t)}$$ is a standard Brownian motion under the measure $\PP^0$. \\
\indent Following the same calculations in \cite{JP2018}, we are led to the following reformulation of the optimal stopping problem in \eqref{osp}
\bel{value}\bar V(x,\psi)=\inf_{\tau\in {\bf S}} \EE_{x,\psi}^0\bracks{1\wedge \varPhi_t+c\int_0^\tau ( 1+\varPhi_t)dt}.\eel
We define the continuation region $C$ and stopping region $D$ by
$$C:=\Big\{(x,\varphi)\in\cX\times(0,\infty):\bar V(x,\varphi)< 1\wedge \varphi\Big\}$$
and $$D:=\Big\{(x,\varphi)\in\cX\times(0,\infty):\bar V(x,\varphi)= 1\wedge \varphi\Big\}.$$
We shall denote by $S$ the boundary between $C$ and $D$. We are now prepared to state this section's key theorem.

\begin{theorem}\label{ospt}
Suppose that either Assumption \ref{time-changeok} or Assumption \ref{time-changeok-2} holds.

{\rm (I)} For any initial couple $(x,\psi)$,  the stopping time  $$\tau^{x,\psi}_D=\inf\{t\geq 0:\psi L_t\notin (l^*_0(X_t),l^*_1(X_t))\},$$ is finite almost surely and is optimal for the the optimal stopping problem in \eqref{osp},  where \bel{l0l1}\left\{\ba{ll}\ad l^*_0(x):=\sup\Big\{\psi\leq 1:  \bar V(x,\psi)=\psi\Big\},\\
\ns\ad l^*_1(x):=\inf\Big\{\psi\geq 1:  \bar V(x,\psi)=1\Big\}.  \ea\right.\eel
Further, we have that $l_0^*(\cdot)<1<l_1^*(\cdot)$.\smallskip

{\rm (II)} Define $A(\varphi_0):=\{(x,\varphi):(x,\varphi/\varphi_0)\in A\}$ for $A=D,S$.  Let $A^\circ(\varphi_0)$ denote the interior of $A(\varphi_0)$.
The  boundary of $D(\varphi_0)$ is probabilistic regular (cf. \cite[p.245]{KA1998})  in the sense that 
$$\PP^\varphi_{x}(\t_{ D^\circ(\varphi_0)}=0)=1\text{ for any }(x,\varphi)\in S(\varphi_0),$$
where $\t_{ D^\circ(\varphi_0)}$ is the entry time of $(X,L)$ to  $D^\circ(\varphi_0)$.\smallskip

{\rm (III)} Given $x_0\in\cX$,  for any $\varphi_0\in (l_0^*(x_0),l_1^*(x_0))$, let $$\tau^*(\varphi_0):=\inf\{t\geq 0:\varphi_0 L_t\notin (l^*_0(X_t),l^*_1(X_t))\}.$$
It follows that
$$\EE^0_{x_0}\int_0^{\tau^*(\varphi_0)}(1+L_t)dt<\infty.$$

\end{theorem}

So as not to detract from the flow of the manuscript, we postpone the proof of Theorem \ref{ospt} to Appendix \ref{AA}. Some remarks, however, are now in order. 

\begin{remark} \rm
 Statement (II) derives the probabilistic regularity of the stopping boundaries for the optimal stopping problem at hand. Probabilistic regularity is not proved by \cite{Gape2011}. As we shall see in Section \ref{sec:mst}, probabilistic regularity of the stopping boundaries plays a key role in solving the minimax Wiener sequential testing problem in \eqref{minimaxV1}. Statement (III) follows from Statement (II) and guarantees the finiteness of  $\bar J$.
\end{remark}

The optimal stopping boundaries given by Theorem \ref{ospt} are depicted in Figure \ref{Fig.main1} below.

\begin{figure}[ht!]
	\centering
	\includegraphics[width=0.5\textwidth]{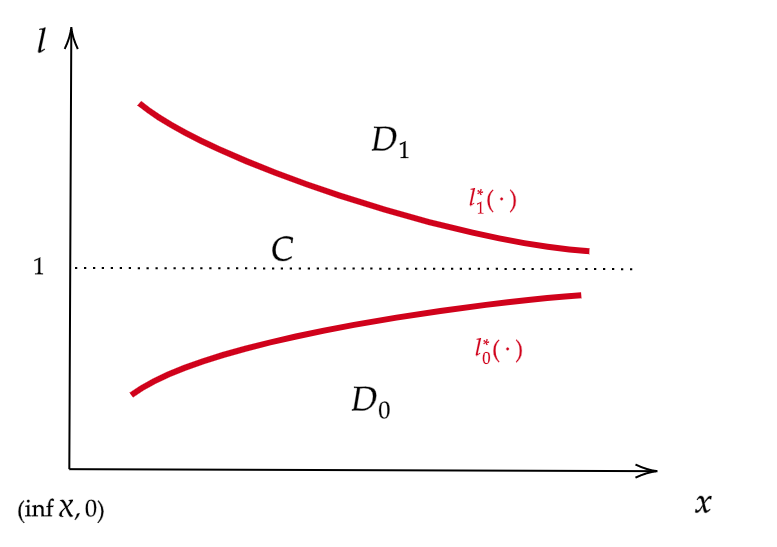}
	\caption{The optimal stopping boundary  under Assumption \ref{time-changeok}.} 
	\label{Fig.main1} 
\end{figure}

\section{Minimax sequential testing}\label{sec:mst}
This section is devoted to solving the minimax Wiener sequential testing problem in \eqref{minimaxV1}. In the sequel, let $x_0$ be any fixed initial state in $\cX$. In order to solve the optimal stopping problem in \eqref{minimaxV1}, we shall find the \textit{least favorable distribution}  $\varphi_0$ in the sense that, for any $\psi\geq 0$,
\bel{lst}\bar J(x_0,\psi;\tau^*(\varphi_0))\leq \bar J\parens{x_0,\varphi_0;\tau^*(\varphi_0)}.\eel
We now introduce the key theorem of this paper, which reveals a characterization of the least favorable distribution for the minimax Wiener sequential testing problem.
\begin{theorem}\label{thm41} Suppose either Assumption \ref{time-changeok} or \ref{time-changeok-2} holds.  Given $x_0\in \cX$,
the set of all least favorable distributions for the optimal stopping problem in \eqref{minimaxV1} is 
\bel{defPsi0} \Psi_0:=\Big\{\varphi\in(l_0^*(x_0),l_1^*(x_0)):\frac{\partial \bar J}{\partial \psi }(x_0,\varphi;\tau^*(\varphi))=0\Big\},\eel
where  $l_0^*(\cdot)$ and $l_1^*(\cdot)$ are defined in \eqref{l0l1}. Further, $\Psi_0\neq\emptyset$. For any $\varphi_0\in\Psi_0$, the corresponding optimal stopping time is  $\tau^*(\varphi_0)$, which has been given in Statement (III) of Theorem \ref{ospt}.
\end{theorem}

\begin{proof}

We begin by invoking Statement (III) in Theorem \ref{ospt}. This tells us that, for all  $\varphi_0\in(l_0^*(x_0),l_1^*(x_0)),$
$$\bar J(x_0,\cdot;\tau^*(\varphi_0))<\infty,$$
and thus both sides of the inequality in \eqref{lst} are finite.\smallskip

Employing the definition of $\bar J$ in \eqref{barJ}, we write
$$\bar J\parens{x_0,\psi;\tau^*(\varphi_0)}=\frac1{1+\psi}\parens{\EE_{x_0}^{0}\bracks{1\wedge (\psi L_{\tau^*(\varphi_0)})}+c\EE_{x_0}^0\int_0^{\tau^*(\varphi_0)} \parens{1+\psi{L_t}}dt}.$$
Note that 
$$\ba{ll}\ad\lim_{\delta\rightarrow 0^+}\frac1\delta\bracks{\EE_{x_0}^{0}\bracks{1\wedge \parens{\parens{\psi+\delta} L_{\tau^*(\varphi_0)}}}-\EE_{x_0}^{0}\bracks{1\wedge \parens{\psi L_{\tau^*(\varphi_0)}}}}\\
        \ns\ad=\lim_{\delta\rightarrow 0^+}\frac1\delta\bracks{\delta\EE_{x_0}^{0}\parens{L_{\tau^*(\varphi_0)}I\bracks{L_{\tau^*(\varphi_0)}\leq \frac1{\psi+\delta}}}+\EE_{x_0}^{0}\parens{\parens{1-\psi L_{\tau^*(\varphi_0)}}I\bracks{\frac1{\psi+\delta}< L_{\tau^*(\varphi_0)}<\frac1\psi}}}\\
\ns\ad=\EE_{x_0}^{0}\parens{L_{\tau^*(\varphi_0)}I\bracks{L_{\tau^*(\varphi_0)}< \frac1{\psi}}}.\ea$$
Taking the right derivative of $\bar J$ with respect to $\psi$, we obtain 
$$\ba{ll}\ad \frac{\partial^+ \bar J}{\partial \psi }\parens{x_0,\psi;\tau^*(\varphi_0)}\\
\ns\ad=\frac1{(1+\psi)^2}\parens{(1+\psi)\EE_{x_0}^0\bracks{L_{\tau^*(\varphi_0)}I\bracks{\psi L_{\tau^*(\varphi_0)}<1 }}-\EE_{x_0}^{0}\bracks{1\wedge \parens{\psi L_{\tau^*(\varphi_0)}}}-c\EE_{x_0}^0\int_0^{\tau^*(\varphi_0)}(1-L_t) dt}\\
\ns\ad=\frac1{(1+\psi)^2}\(\EE_{x_0}^0\bracks{L_{\tau^*(\varphi_0)}I\bracks{\psi L_{\tau^*(\varphi_0)}<1)}}-\PP_{x_0}^0(\psi L_{\tau^*(\varphi_0)}\geq 1)-c\EE_{x_0}^0\int_0^{\tau^*(\varphi_0)}(1-L_t)dt\). \ea$$
Note that 
$$\ba{ll}\ad\lim_{\delta\rightarrow 0^+}\frac1\delta\bracks{\EE_{x_0}^{0}\bracks{1\wedge (\psi L_{\tau^*(\varphi_0)})}-\EE_{x_0}^{0}[1\wedge ((\psi-\delta) L_{\tau^*(\varphi_0)})]}\\
\ns\ad=\lim_{\delta\rightarrow 0^+}\frac1\delta\bracks{\delta\EE_{x_0}^{0}\parens{L_{\tau^*(\varphi_0)}I\bracks{ L_{\tau^*(\varphi_0)}\leq \frac1{\psi}}}+\EE_{x_0}^{0}\parens{(1-\psi L_{\tau^*{\varphi_0}})I\bracks{\frac1{\psi}< L_{\tau^*(\varphi_0)}<\frac1{\psi-\delta}}}}\\
\ns\ad=\EE_{x_0}^{0}\parens{L_{\tau^*(\varphi_0)}I\bracks{ L_{\tau^*(\varphi_0)}\leq  \frac1{\psi}}}.\ea$$
Taking  the left derivative of $\bar J$ with respect to $\psi$, we obtain
$$\ba{ll}\ad \frac{\partial^- \bar J}{\partial \psi }(x_0,\psi;\tau^*(\varphi_0))\\
\ns\ad=\frac1{(1+\psi)^2}\parens{(1+\psi)\EE_{x_0}^0\parens{L_{\tau^*(\varphi_0)}I\bracks{\psi L_{\tau^*(\varphi_0)}\leq 1 }}-\EE_{x_0}^{0}\bracks{1\wedge \parens{\psi L_{\tau^*(\varphi_0)}}}-c\EE_{x_0}^0\int_0^{\tau^*(\varphi_0)}(1-L_t) dt}\\
\ns\ad=\frac1{(1+\psi)^2}\parens{\EE_{x_0}^0\bracks{L_{\tau^*(\varphi_0)}I\bracks{\psi L_{\tau^*(\varphi_0)}\leq 1}} -\PP_{x_0}^0\parens{\psi L_{\tau^*(\varphi_0)}> 1}-c\EE_{x_0}^0\int_0^{\tau^*(\varphi_0)}\parens{1-L_t}dt}.\ea$$
One may easily see that the partial derivatives $\frac{\partial^\pm \bar J}{\partial \psi }\parens{x_0,\cdot;\tau^*(\varphi_0)}$ are decreasing with respect to $\psi$. Therefore  $\varphi_0$ is a maximum point in \eqref{lst} if and only if the following inequality is satisfied
$$ \frac{\partial^+ \bar J}{\partial \psi }\parens{x_0,\varphi_0;\tau^*(\varphi_0)} \leq 0 \leq \frac{\partial^- \bar J}{\partial \psi }(x_0,\varphi_0;\tau^*(\varphi_0)).$$
 
\indent We proceed to prove the existence of this $\varphi_0$. 
Recall that $(X_t,L_t)$ has continuous paths and that, from Statement (I) of Theorem \ref{ospt}, $$l^*_0(\cdot)<1<l^*_1(\cdot).$$  Then, at $t=\tau^*(\varphi_0)$, we have that $$1\neq \varphi_0L_{\tau^*(\varphi_0)}=l^*_0(X_{\tau^*(\varphi_0)})\,\, \text{or}\,\, l^*_1(X_{\tau^*(\varphi_0)}).$$
We may then write
\bel{nomass}\frac{\partial^- \bar J}{\partial \psi }(x_0,\varphi_0;\tau^*(\varphi_0))-\frac{\partial^+ \bar J}{\partial \psi }(x_0,\varphi_0;\tau^*(\varphi_0))=(1+1/{\varphi_0})\PP^0\parens{L_{\tau^*(\varphi_0)}=1/\varphi_0}=0.\eel
We then have, for any $\varphi_0\in (l_0^*(x_0),l_1^*(x_0))$,
$$\frac{\partial \bar J}{\partial \psi }(x_0,\varphi_0;\tau^*(\varphi_0))=\frac{\partial^\pm \bar J}{\partial \psi }(x_0,\varphi_0;\tau^*(\varphi_0)).$$
\indent We consider two possible cases for the values of $\frac{\partial \bar J}{\partial \psi }(x_0,\varphi_0;\tau^*(\varphi_0))$. If $\varphi_0=l_1^*(x_0)$, then $\tau^*(l_1^*(x_0))=0$ and $L_{\tau^*(l_1^*(x_0))}=1$. We then have $$\frac{\partial \bar J}{\partial \psi }\parens{x_0,\varphi_0;\tau^*(\varphi_0)}=-1.$$ However, if $\varphi_0=l_0^*(x_0)$, then $\tau^*(l_0^*(x_0))=0$ and  $L_{\tau^*(l_1^*(x_0))}=1$. We then have $$\frac{\partial \bar J}{\partial \psi }(x_0,\varphi_0;\tau^*(\varphi_0))=1.$$ 
\indent It now suffices to prove $\frac{\partial \bar J}{\partial \psi }(x_0,\varphi_0;\tau^*(\varphi_0))$ is continuous with respect to $\varphi_0$ in  $[l_0^*(x_0),l_1^*(x_0)]$. Let $\varphi_0$ be fixed n  $[l_0^*(x_0),l_1^*(x_0)]$ and define
\begin{align*}&\overline\th(\varphi):=\inf\Big\{t\geq 0: (X_t, L_t)\notin C(\varphi)\cup C(\varphi_0) \Big\},\\& \underline\th(\varphi):=\inf\Big\{t\geq 0: (X_t, L_t)\notin C(\varphi)\cap C(\varphi_0) \Big\}.\end{align*}
Taking $\varphi\rightarrow\varphi_0$, it immediately follows that $$ \big[ C(\varphi)\cap C(\varphi_0)\big]\subset C(\varphi_0)\text{ with } \big[ C(\varphi)\cap
C(\varphi_0)\big]\uparrow C(\varphi_0),$$ and that  $$ \big[C(\varphi)\cup C(\varphi_0)\big]\supset C(\varphi_0)\text{ with } \big[C(\varphi)\cup C(\varphi_0)\big]\downarrow  C(\varphi_0).$$ Since $\tau^*(\varphi_0)<\infty$ almost surely, $\underline\th(\varphi)\uparrow \tau^*(\varphi_0)$ almost surely. Invoking the probabilistic regularity of $S(\varphi_0)$ from Statement (II) in Theorem \ref{ospt} yields $$\overline\th(\varphi) \downarrow \{t\geq 0: (X_t, L_t)\in  D^\circ(\varphi_0)\}= \tau^*(\varphi_0)\,\text{ a.s.}$$   Therefore when $\varphi\rightarrow\varphi_0$, we have that $\tau^*(\varphi)\rightarrow\tau^*(\varphi_0)$ almost surely. The continuity of \begin{align*} &\frac{\partial \bar J}{\partial \psi }(x_0,\varphi_0;\tau^*(\varphi_0))\\
&=\frac1{(1+\varphi_0)^2}\(\EE_{x_0}^0\bracks{L_{\tau^*(\varphi_0)}I(\varphi_0 L_{\tau^*(\varphi_0)}\leq 1)}-\PP_{x_0}^0(\varphi_0 L_{\tau^*(\varphi_0)}> 1)-c\EE_{x_0}^0\int_0^{\tau^*(\varphi_0)}(1-L_t)dt\),\end{align*} with respect to $\varphi_0\in [l_0^*(x_0),l_1^*(x_0)]$ holds by recalling the fact that $\varphi_0L_{\tau^*(\varphi_0)}$ does not admit a positive probability at $1$ (as noted in \eqref{nomass}) and by invoking dominated convergence. This completes the proof of existence of $\varphi_0$ when $$ \frac{\partial \bar J}{\partial \psi }(x_0,\varphi_0;\tau^*(\varphi_0))=0.$$
We may now conclude that $\varphi_0$ is a least favorable distribution if and only if $\varphi_0\in\Psi_0$. Moreover, $\Psi_0\neq \emptyset$. This completes the proof.
\end{proof}

\begin{remark} \rm We pause to compare the optimal stopping decisions for the optimal stopping problem in \eqref{osp} and for the optimal stopping problem in \eqref{minimaxV1}. In the setting of the optimal stopping problem in \eqref{osp}, an initial priori distribution $\psi$ is given and $\tau^*(\psi)$ is the optimal stopping time as presented in Theorem \ref{ospt}. However, in the setting of the minimax Wiener sequential testing problem given in \eqref{minimaxV1}, the priori distribution $\psi$ is not given. This requires us to find a least favorable distribution $\psi=\varphi_0\in\Psi_0$. Once we find this least variable distribution, $\tau^*(\varphi_0)$ then becomes the optimal stopping time for the optimal stopping problem in \eqref{minimaxV1}.
\end{remark}

\begin{remark} \rm
Although it seems rather impossible to obtain an analytical closed form for $(\varphi_0,l_0^*(\cdot),l_1^*(\cdot))$, one may of course easily turn to numerics. To find the value function $V$ for the optimal stopping problem in \eqref{osp}, one would solve the following Hamilton-Jacobi equation 
$$\min(\LL^0 V(x,\psi)+c(1+\psi),  1\wedge \psi -V(x,\psi))=0,$$
and proceed to find the optimal stopping boundary $(l_0^*(\cdot),l_1^*(\cdot))$ through \eqref{l0l1}. A bisection search may then be employed to numerically find  $\varphi_0\in \Psi_0$.
\end{remark}

\subsection{Examples}\label{sec:exp}

In this section, we provide two concrete examples to illustrate Theorem \ref{thm41}.\\

\noindent \textbf{Example 1}.\\ We revisit the work of Johnson and Peskir \cite{JP2018}. Consider the following $\alpha$-dimensional Bessel process in $\cX=(0,\infty)$
$$dX_t=\frac{\alpha-1}{X_t}dt+dB_t \text{ with } X_0=x_0>0.$$
The dimension $\alpha$ may be either $\delta_0$ and $\delta_1$. Without loss of generality, we may assume that $\delta_1>\delta_0$. In this case: $\mu_0(x)=\frac{\delta_0 -1}{2x}$, $\mu_1(x)=\frac{\delta_1 -1}{2x}$ and $\sigma(x)=1$. We then calculate $\rho$ as
$$\rho^2(x)=\frac{(\delta_1-\delta_0)^2}{4x^2},$$
which is clearly a decreasing function.
The function $K(x)$ in \eqref{DefK} satisfies
$$K(x)=\frac {\delta_0-2}{\delta_1-\delta_0}>-\frac12\text{ if }\delta_1+\delta_0>4.$$
Given any $x_0>0$, if $\delta_0+\delta_1>4$, we may invoke Theorem \ref{thm41} to conclude that the least favorable distribution for the minimax Wiener sequential testing problem in \eqref{minimaxV1} exists and is given by \eqref{defPsi0}.\\

\noindent \textbf{Example 2}. \\
We revisit the work of Gapeev and Shiryaev (\cite{Gape2011}). Let
\begin{align*}
&\mu_i(x)=\frac{\eta_i \sigma^2(x)}{x}, 
\end{align*}
for some $\eta_1,\eta_2\in\dbR$ and for some  $\sigma(\cd)$ such that 
$\cX=(0,\infty).$
In this case, 
$$K(x)+\frac12=\frac{\eta_0+\eta_1-1}{2(\eta_1-\eta_0)}.$$
Theorem \ref{thm41} yields that if either
\begin{align}\label{choice12}
&\frac{\eta_0+\eta_1-1}{2(\eta_1-\eta_0)}<0,\q \frac{\sigma^2(x)}{x^2}\text{ is increasing on $(0,\infty)$},\\
\text{ or }\,\,&\frac{\eta_0+\eta_1-1}{2(\eta_1-\eta_0)}>0,\q \frac{\sigma^2(x)}{x^2}\text{ is decreasing on $(0,\infty)$,}\nonumber \end{align}
and if
\begin{align}\label{distingsui}\int_0^\cdot\frac{\sigma^2(X_t)}{X_t^2}dt\text { is strictly increasing to $\infty$ under either $\PP^0$ or $\PP^\infty$, }
\end{align}
 then the minimax sequential testing problem admits a  least favorable distribution satisfying \eqref{defPsi0} (our measure $\PP^0$ is the same as the measure $\PP_0$ in  \cite{Gape2011} and our measure $\PP^\infty$ is the same as the measure $\PP_1$ in \cite{Gape2011}). The choice in \eqref{choice12} only depends on the monotonicity of $\sigma^2(x)/x^2$. For example, consider the case where $\sigma^2(x)/{x^2}$ is increasing on $(0,\infty)$. If $\eta_0+\eta_1<1$, we would choose $\eta_1>\eta_0$ when specifying the hypotheses $H_0$ and $H_1$. Similarly, if $\eta_0+\eta_1 \geq 1$, we would choose $\eta_0<\eta_1$ when specifying the hypotheses $H_0$ and $H_1$. This means that  if $\sigma^2(x)/x^2$ is monotone on $(0,\infty)$, \eqref{choice12} holds. Should \eqref{distingsui} also hold, then there exists a least favorable distribution which satisfies \eqref{defPsi0}.

\section{The case of constant SNR }\label{sec:SNRC}

 In all of the previous sections, we have assumed that the signal-to-noise ratio $$ \rho(x)=\frac{\mu_1(x)-\mu_0(x)}{\sigma(x)}$$ is non-constant. We will now assume, for all $x\in\cX$, that the signal-to-noise ratio is constant and equal to $\rho_0$. We proceed to consider the minimax Wiener sequential testing problem in \eqref{minimaxV1} under the assumption of constant SNR. Trivially, Assumption \ref{DefK} and Assumption \ref{DefK2} from Section \ref{sec:osp} are now irrelevant.\\
\indent In order to derive the form of the least favorable distribution, we will consider  a generalized performance functional $U$ defined as
$$U(x,\psi;\tau,d_\tau)=\PP_x^\psi(d_\tau\neq \th)+\EE^\psi_x\int_0^{\tau}f(L_t)dt,$$
where $f:(0,\infty)\mapsto \dbR^+$ is a smooth running cost function of the likelihood ratio.
When $f=1$, the performance functional $U$ reduces to the performance functional $J$ in \eqref{per-1}.\\
\indent In the setting of constant SNR, the likelihood ratio process $L_t$ (as defined in \eqref{liklihood2}) is a strong Markov process. The optimal decision for the optimal stopping problem can then be determined entirely by the likelihood ratio process $L_t$.
In order to minimize the probability of false detection, we select the decision $\bar d$ as in \eqref{decisionprocedure}. Moreover, note that $U$, $l_1^*(x)$ and $l_0^*(x)$ are all independent of $x$. This leads us to formulate the following performance functional independent of   $x$ and $d$ $$\bar U(\psi;\tau)=U(x,\psi;\tau,\bar d).$$
The minimax Wiener sequential testing problem 
in this case then reduces to the following optimal stopping problem. \smallskip

\noindent{\bf OSP 3}: Given $x_0\in\cX$, find a $\tau^*\in{\bf S} $ such that \bel{minimaxSNRC}\sup_{\psi\geq 0}\bar U(\psi;\tau^*)=\inf_{\tau\in{\bf S}}\sup_{\psi\geq 0}\bar U(\psi;\tau).\eel

We proceed to solve the optimal stopping problem in \eqref{minimaxSNRC}. We break up the solution into two intermediate steps.\smallskip

(1) {\it Optimal strategy.} Suppose, for some $\varphi_0>0$, that the initial value of $\psi$ is $\varphi_0$. We would then seek to find  the optimal stopping time $\tau^*$ for the following optimal stopping problem  
\begin{equation}
\bar U(\varphi_0;\tau^*)=\inf_{\tau\in {\bf S}}\bar U(\varphi_0;\tau).
\end{equation}
The candidate optimal stopping time to be verified is the stopping time $\tau^*(\varphi_0)$.
Using standard tools from the theory of optimal stopping for diffusions (see \cite{PS}), the form of $\tau^*$ can be found from the following free-boundary problem 
\bel{fbp}\left\{\ba{ll}\ad \LL^0W+(1+l)f=0 \text{ for }l^*_0< l< l^*_1;\\
\ns\ad W(l^*_0)=l^*_0, W(l^*_1)=1;\\
\ns\ad W'(l^*_0)=1, W'(l^*_1)=0,\ea\right.\eel
where $$\LL^0 W(l):=-\frac {\rho_0^2}{2}l^2W''(l).$$
\indent We proceed by defining the function $M(\cdot)$
$$M(l)=-\frac{2}{\rho_0^2}\int_1^l\int_1^v \frac{f^2(u)}{ u^2}dudv,$$
which is the solution to 
$$\LL^0 M= (1+l) f\text{ with }M(1)=M'(1)=0. $$
It is straightforward to see that if $l^*_0< l< l^*_1$, we have, for some $A,B$,
$$W(l)=M(l)+A l+B.$$ The boundary conditions then become
$$\left\{\ba{ll}\ad  M(l^*_0)+A l^*_0+B=l^*_0, ~M(l^*_1)+A l^*_1+B=1,\\
\ns\ad M'(l^*_0)+A=1,~M'(l^*_1)+A=0,\ea\right.$$
from which the explicit value of $(l^*_0, l^*_1)$ can be identified for given $\rho_0$ and  $f$. If $W(l)\leq 1\wedge l$ on $(l_0^*,l_1^*)$, we may then apply It\^o's formula. The stopping time $\tau^*$  defined by
$$\tau^*=\inf\{t\geq 0: L_t\notin (l^*_0/\varphi_0,l^*_1/\varphi_0)\}=\tau^*(\varphi_0),$$
will then be optimal.\smallskip

(2) {\it Verification of least favorable distribution}. To verify that $\tau^*(\varphi_0)$ is optimal for the optimal stopping problem in \eqref{minimaxSNRC}, it is  equivalent to check that, for all $\psi>0$, the following inequality holds
$$\bar U(\psi;\tau^*(\varphi_0))\leq  \bar U(\varphi_0;\tau^*(\varphi_0)).$$

Lengthy but straightforward calculations yield 
$$\bar  U(\psi;\tau^*(\varphi_0))=\left\{\ba{ll} \ad \frac{1\wedge \psi }{1+\psi}+\frac 1{1+\psi}\parens{\frac{\varphi_0-l^*_0}{l^*_1-l^*_0}M(l^*_1/\varphi_0)+\frac{l^*_1-\varphi_0}{l^*_1-l^*_0}M(l^*_0/\varphi_0)}\\
\ns\ad\q+\frac \psi {1+\psi}\parens{\frac{l^*_1(\varphi_0-l^*_0)}{\varphi_0(l^*_1-l^*_0)}M(\varphi_0/l^*_1)+\frac{l^*_0(l^*_1-\varphi_0)}{\varphi_0(l^*_1-l^*_0)}M(\varphi_0/l^*_0)},\\[2mm]
\ns\ad\qq\qq\q\text{ if } \varphi_0/\psi\leq l^*_0 \text{ or }\varphi_0/\psi\geq l^*_1;\\[2mm]
\ns\ad \frac 1{1+\psi}\parens{\frac{\varphi_0-l^*_0}{l^*_1-l^*_0}+\frac{\varphi_0-l^*_0}{l^*_1-l^*_0}M(l^*_1/\varphi_0)+\frac{l^*_1-\varphi_0}{l^*_1-l^*_0}M(l^*_0/\varphi_0)}\\
\ns\ad\q+\frac \psi {1+\psi}\parens{\frac{l^*_0(l^*_1-\varphi_0)}{\varphi_0(l^*_1-l^*_0)}+\frac{l^*_1(\varphi_0-l^*_0)}{\varphi_0(l^*_1-l^*_0)}M(\varphi_0/l^*_1)+\frac{l^*_0(l^*_1-\varphi_0)}{\varphi_0(l^*_1-l^*_0)}M(\varphi_0/l^*_0)} \\[2mm]
\ns\ad\qq\qq\qq\text{ if }\varphi_0/l^*_1<\psi<\varphi_0/l^*_0.\ea\right.$$
We then have that $\psi=\varphi_0$ is the maximum point of $\bar U(\psi;\tau_0^*(\varphi_0))$ if and only if
\bel{case1}\ba{ll}\ns\ad\frac{l^*_0(l^*_1-\varphi_0)}{\varphi_0(l^*_1-l^*_0)}+\frac{l^*_1(\varphi_0-l^*_0)}{\varphi_0(l^*_1-l^*_0)}M(\varphi_0/l^*_1)+\frac{l^*_0(l^*_1-\varphi_0)}{\varphi_0(l^*_1-l^*_0)}M(\varphi_0/l^*_0)\\
\ns\ad=  \frac{\varphi_0-l^*_0}{l^*_1-l^*_0}+\frac{\varphi_0-l^*_0}{l^*_1-l^*_0}M(l^*_1/\varphi_0)+\frac{l^*_1-\varphi_0}{l^*_1-l^*_0}M(l^*_0/\varphi_0).\ea\eel
Note that, for $\varphi_0=l_0^*$,
$${\rm LHS\,\, of \,\, \eqref{case1}}=1+M(1)>M(1)={\rm RHS\,\,of\,\, \eqref{case1}}.$$
For $\varphi_0=l_1^*$,
$${\rm LHS\,\, of \,\, \eqref{case1}}=M(1)<1+M(1)={\rm RHS \,\, of \,\, \eqref{case1}}.$$
By the mean-value theorem, there exists a $\varphi_0\in(l^*_0,l^*_1) $ such that \eqref{case1} holds. Theorem \ref{thm51} below then immediately follows.

\begin{theorem} \label{thm51} Suppose the following free-boundary problem
$$\left\{\ba{ll}\ad\LL^0 M= (1+l) f\text{ for }l\in (l_0^*,l_1^*),\\
\ns\ad  M(l^*_0)+A l^*_0+B=l^*_0, ~M(l^*_1)+A l^*_1+B=1,\\
\ns\ad M'(l^*_0)+A=1,~M'(l^*_1)+A=0,\\
\ns\ad M(l)+A l+B\leq 1\wedge l \text{ for }l\in (l_0^*,l_1^*),\\
&0<l_0^*<1<l_1^*,\ea\right.$$
 admits a solution $(M, A,B,l_0^*, l_1^*)$. Then $\varphi_0\in(l_0^*,l_1^*)$  solving \eqref{case1} is a least favorable distribution and  the stopping time $\tau_0^*=\inf\{t\geq 0:\varphi_0 L_t\notin (l_0^*,l_1^*)\}$ is optimal for the optimal stopping problem in \eqref{minimaxSNRC}.
\end{theorem}

\subsection{The Case where both SNR and $f$ are constant}
We now examine the special case where both SNR and $f$ 
are constant. In this case, we can explicitly find a closed form solution for the least favorable distribution. 

\begin{proposition}[Symmetric case when SNR is a constant]
Suppose the following free-boundary problem
$$\left\{\ba{ll}\ad\LL^0 M= (1+l) f\text{ for }l\in (l_0^*,l_1^*),\\\ns\ad  M(l^*_0)+A l^*_0+B=l^*_0, ~M(l^*_1)+A l^*_1+B=1,\\
\ns\ad M'(l^*_0)+A=1,~M'(l^*_1)+A=0,\\
\ns\ad M(l)+A l+B\leq 1\wedge l \text{ for }l\in (l_0^*,l_1^*),\\
&0<l_0^*<1<l_1^*,\ea\right.$$
 admits a solution $(M, A,B,l_0^*, l_1^*)$. If  $f(l)=f(1/l)$ for any $l\in(0,\infty)$, then $\varphi_0=1$ is a least favorable distribution and the first exit time of $L_t$ from $(l^*_0,l^*_1)$ is optimal for the optimal stopping problem in \eqref{minimaxSNRC}. When $f$ is constant, the free boundary problem admits a solution and $\varphi_0=1$ is the least favorable distribution.
\end{proposition}

\begin{proof}
We begin by verifying that \eqref{case1} holds for $\varphi_0=1$ in the symmetric case. We claim that 
$l^*_0=1/l_1^*.$ Recall that $\tau_0^*=\inf\{t\geq 0: \psi L_t\notin (l_0^*,l_1^*)\}$ is the solution to the optimal stopping problem
$$R(\psi):=\inf_{\tau\in{\bf S}}\bar U(\psi;\tau).$$
We have 
\bel{Rvalue} R(l_0^*)=l_0^* \text{ and }R(l_1^*)=1.\eel
Since the probability law of the likelihood ratio process $L$ under $\PP^0$ is same as the probability law of $L^{-1}$ under the measure $\PP^\infty$, 
$$\ba{ll}\ad R(\psi)=\inf_{\tau\in{\bf S}}\frac{\psi}{1+\psi}\EE^0\[ 1\wedge (\psi L_{\tau})+c\int_0^\tau (1+\psi L_t) f(L_t)dt\]\\
\ns\ad=\inf_{\tau\in{\bf S}}\frac{\psi}{1+\psi}\EE^\infty\[ 1\wedge (\psi L^{-1}_{\tau})+c\int_0^\tau (1+\psi L^{-1}_t) f(L_t)dt\]\\
\ns\ad=\psi\inf_{\tau\in{\bf S}}\frac{\psi^{-1}}{1+\psi^{-1}}\EE^\infty\[ 1\wedge (\psi L^{-1}_{\tau})+c\int_0^\tau (1+\psi L^{-1}_t) f(L^{-1}_t)dt\]\\
\ns\ad=\psi R(\psi^{-1}).\ea$$
We then have that $R(1/l_0^*)=1=R(l_1^*)$, which indicates that $l_0^*=1/l_1^*$. Plugging the equality $l_0^*=1/l_1^*$ into \eqref{case1} with $\varphi_0=1$, we see that
$$\text{LHS of \eqref{case1}}=\frac{l_0^*}{1+l_0^*}+\frac{1}{1+l_0^*}M(l^*_0)+\frac{l_0^*}{1+l_0^*}M(1/l_0^*)=\text{RHS of \eqref{case1}}.$$
This says that $\varphi_0=1$ is the least favorable distribution.\\
\indent When $f$ is constant, Peskir and Shiryaev \cite[p. 290]{PS} prove that the free boundary problem has a unique solution. We may thus conclude that $\varphi_0=1$ is the least favorable distribution in the case of constant SNR and constant $f$. This completes the proof. 
\end{proof}
\bigskip \bigskip

\noindent \textbf{Acknowledgements} Philip A. Ernst gratefully acknowledges support from the Army Research Office (ARO-YIP-71636-MA), the National Science Foundation (DMS-2311306), the Office of Naval Research (N00014-21-1-2672), and the Royal Society Wolfson Fellowship (RSWF$\backslash$R2$\backslash$222005). Hongwei Mei gratefully acknowledges support the Simons Foundation's Travel Support for Mathematicians Program (No. 00002835).



\newpage

\newpage

\begin{center}
	
\end{center}

\newpage 

\appendix
\section{Appendix: Proof of Theorem \ref{ospt}}\label{AA}

The purpose of this appendix is to prove Theorem \ref{ospt} from Section \ref{sec:osp}. We need only prove the theorem under  Assumption \ref{time-changeok} as the proof of Theorem \ref{ospt} under Assumption \ref{time-changeok-2} is completely symmetric. For the sake of convenience, in the sequel, we shall omit the superscript 0 in $\EE^0$ or $\PP^0$.\smallskip

\indent We begin with the proof of Statement (I) of Theorem \ref{ospt}. For ease of exposition, we shall split the proof into several intermediate steps.
\smallskip

 (1) {\it Time-change}. 
 Let us define $T_t\geq 0$ by
$$\int_0^{T_t}\rho^2(X_s)ds=t.$$ Note that $T_t$ is strictly increasing and is uniquely defined. Let 
$\hat X_t:=X_{T_t}$ and $\hat \varPhi_t:=\varPhi_{T_t}.$
It is straightforward to see that 
${dT_t}/{dt}=\rho^{-2}(\hat X_t).$
Employing It\^o's formula, we obtain
\bel{hatSDE}\left\{\ba{ll}\ad d\hat X_t=\mu_0(\hat X_t)\rho^{-2}(\hat X_t)dt+\sigma(\hat X_t)\rho^{-1}(\hat X_t)d\hat B_t\\
\ns\ad  d\hat \varPhi_t=\hat\varPhi_t d\hat B_t,\ea\right.\eel
where $\hat B_t=\int_0^{T_t}\rho(X_t)d\wt B_t$ is a standard Brownian motion under the measure $\PP$ by its L\'evy characterization.
Then the optimal stopping problem in \eqref{value} is equivalent to
\bel{hatV}\ba{ll}  V(x,\varphi)\ad=\inf_{\t\in\cT_{\hat X,\hat\varPhi}}\EE_{x,\varphi}\[\int_0^\tau H(\hat X_t,\hat\varPhi_t)dt+G(\hat\varPhi_\t)\],\ea\eel
where $H(x,\varphi):=c\rho^{-2}(x)(1+\varphi)$, $G(\varphi):=(1\wedge \varphi)$, $\cT_{\hat X,\hat\varPhi}$ is the set of stopping times of $(\hat X,\hat\varPhi)$, and the subscript $x,\varphi$ under $\EE$ stands for the initial value of $\widehat X_0,\widehat\varPhi_0$.  Since $\hat\varPhi$ is an exponential martingale,  $\hat\varPhi_t\rightarrow0$ almost surely as $t\rightarrow\infty$. For simplicity, we write $\hat\mu=\mu_0\rho^{-2}$ and $\hat\sigma=\sigma\rho^{-1}$. The infinitesimal generator of $(\hat X,\hat\varPhi)$ is given by 
$$\hat \LL=\hat \mu(x)\partial_x+\frac12\hat\sigma^2(x)\partial_{xx}+\hat\sigma(x)\varphi\partial_{x\varphi}+\frac12\varphi^2\partial_{\varphi\varphi}.$$
\indent We proceed to study the optimal stopping problem in \eqref{hatV},  which has the same optimal stopping boundary as  the optimal stopping problem in \eqref{value}. We begin with Lemma \ref{Time} below.

\begin{lemma}\label{Time} 
The following three statements hold:\\
	{\rm (1)} The random variable $\log(\varphi^{-1}\sup_{t\geq 0}\hat\varPhi_t)$ is exponentially distributed with mean 1 under $\PP_\varphi$.\smallskip
	
		\noindent {\rm (2)} It follows that \bel{h}h(\varphi):=\EE_\varphi \sup_{t\geq 0}G(\hat\varPhi_t)=\begin{cases}
		1,&\text{ for }\varphi\geq1,\\
		\varphi(1-\log\varphi),&\text{ for }\varphi<1.
		\end{cases}\eel
	{\rm (3)} Let
	$$\gamma_T(\varphi):=\EE_{\varphi}\(\sup_{t\geq T}G(\hat\varPhi_t)\).$$ Then for any $\varphi_0>0$,   
 $\gamma_T(\varphi)$  converges to 0 uniformly in $(0,\varphi_0)$ as $T\rightarrow\infty$.
\end{lemma}

\begin{proof} (1).  Note that $$\log\parens{\varphi^{-1}\sup_{t\geq 0}\hat\varPhi_t}=\sup_{t\geq 0}\(-\frac t2+\hat B_t\).$$
The right-hand side is exponentially distributed with mean 1.\smallskip

(2). Equation \eqref{h}  holds by a straightforward calculation from statement (1).\smallskip

(3).  Note that $\hat\varPhi_t=\varphi\exp\{-t/2+\hat B_t\}$. It follows that 
$$\ba{ll}\ad \EE_{\varphi}\(\sup_{t\geq T}G(\hat\varPhi_t)\)=\EE_\varphi\[\EE\( 1\wedge\sup_{t\geq T}\hat\varPhi_t\Big|\hat\varPhi_T\)\]=\EE_\varphi h(\hat\varPhi_T)
\\
\ns\ad= \PP\(-T/2+\hat B_T+\log\varphi\geq0\)\\
\ns\ad\q+\EE\(\varphi\exp\{-T/2+\hat B_T\}(1+T/2-\hat B_T-\log\varphi)I[-T/2+\hat B_T+\log\varphi<0] \)=\gamma_T(\varphi).\ea$$	
Given any $\varphi>0$, it immediately follows that $\gamma_T(\varphi)\rightarrow0$	 as $T\rightarrow\infty$. Further note that $\gamma_T(\cdot)$ is increasing with respect to $\varphi$. This yields uniform convergence on any finite interval $(0,\varphi_0)$.\end{proof}

(2) {\it Optimal stopping time}. We fix $T>0$ and consider the optimal stopping problem
\bel{hatV-T}\ba{ll}  V^T(t;x,\varphi)\ad=\inf_{0\leq \t\leq T-t}\EE_{x,\varphi}\[\int_0^\tau H(\hat X_t,\hat\varPhi_t)dt+G(\hat\varPhi_\t)\].\ea\eel
By the Feller property of the process $(X,\varPhi)$, it is easy to see that, for all $t\in[0,T]$, $V^T(t;x,\varphi)$ is a continuous function of $(x,\varphi)$. Employing Lemma \ref{Time}, we have
\bel{VTV}\ba{ll}\ad V^T(0;x,\varphi)\geq V^{T+1}(0;x,\varphi)
\geq V(x,\varphi)=\inf_{\t}\EE_{x,\varphi}\[\int_0^\tau H(\hat X_t,\hat\varPhi_t)dt+G(\hat\varPhi_\t)\]\\
\ns\ad\geq\inf_{\t}\EE_{x,\varphi}\[\int_0^{\tau\wedge T} H(\hat X_t,\hat\varPhi_t)dt+G(\hat\varPhi_{\tau\wedge T})+G(\hat\varPhi_{\tau})-G(\hat\varPhi_{\tau\wedge T})\]\\
\ns\ad\geq \inf_{\t}\EE_{x,\varphi}\[\int_0^{\tau\wedge T} H(\hat X_t,\hat\varPhi_t)dt+G(\hat\varPhi_{\tau\wedge T})\]-\EE_\varphi G(\hat\varPhi_T) \geq V^T(0;x,\varphi)-\gamma_T(\varphi).\ea\eel
The above implies that $  \{V^T(0;x,\varphi):T=1,2,\cdots\}$ is a Cauchy sequence in $\sC^1(U)$ for any finite open subset $U\subset\cX\times(0,\infty)$. We now let $T\rightarrow\infty$. By Statement (3) in Lemma \ref{Time}, $V$ is also a continuous function and, for any $\varphi_0<\infty$,
$V^T(0;x,\varphi)$ converges to $V(x,\varphi)$ uniformly on  $\cX\times(0,\varphi_0)$.  \\
\indent We now wish to find the optimal stopping time for $V$ using the sequence of optimal stopping times for $V^T$. To this end, we write
$$C^T(t)=\{(x,\varphi):V^T(t;x,\varphi)< G(\varphi)\}\text{ and }D^T(t)=\{(x,\varphi):V^T(t;x,\varphi)= G(\varphi)\}.$$
We also write 
$$C_\e^T(t)=\{(x,\varphi):V^T(t;x,\varphi)< G(\varphi)-\e\}\text{ and }D_\e^T(t)=\{(x,\varphi):V^T(t;x,\varphi)\geq G(\varphi)-\e\}.$$
Further, let us denote $$\tau_{D_\e^T}:=\inf\{t\in[0,T]:(\hat X_t,\hat\varPhi_t)\in D_\e^T(t) \},$$ with $\inf\emptyset=T$. For any triple $(t;x,\varphi)$ with $t<T$, we have that, for $T'\geq T$,
$$V^T(t;x,\varphi)\geq V^{T'}(t;x,\varphi).$$
This means that $\tau_{D_\e^T}$ increases as  either $T$ increases or $\e$ decreases.
Let $\tau_*=\sup_{\e,T}\tau_{D_\e^T}\leq\tau_{D}$ almost surely. Employing the Snell envelope (see, for example, \cite[Theorem 2.2]{PS}), we have that
$$V^T\parens{t\wedge\tau_{D_\e^T};\hat X_{t\wedge\tau_{D_\e^T}},\hat \varPhi_{t\wedge\tau_{D_\e^T}}}+\int_0^{t\wedge\tau_{D_\e^T}}H(\hat X_t,\hat\varPhi_t)dt,$$
is a martingale. We proceed by writing
\bel{VTVV} \ba{ll}\ad V(x,\varphi)\geq V^T(0;x,\varphi)-\gamma_N(\varphi)\\
\ns\ad=\EE_{x,\varphi}V^T\parens{T;\hat X_{\tau_{D_\e^T}},\hat \varPhi_{\tau_{D_\e^T}}}+\EE_{x,\varphi}\int_0^{\tau_{D_\e^T}}H(\hat X_t,\hat\varPhi_t)dt-\gamma_T(\varphi)\\
\ns\ad= \EE_{x,\varphi}G\parens{\hat \varPhi_{\tau_{D_\e^T}}}+\EE_{x,\varphi}\int_0^{\tau_{D_\e^T}}H(\hat X_t,\hat\varPhi_t)dt-\gamma_T(\varphi),\\ \ea\eel
where we have used the fact that $V^T(T,x,\varphi)=G(\varphi).$ As both $T\rightarrow\infty$ and $\e\rightarrow 0^+$, the last equality in \eqref{VTVV} tends to $$\EE_{x,\varphi}G(\hat \varPhi_{\tau_*})+\EE_{x,\varphi}\int_0^{\tau_*}H(\hat X_t,\hat\varPhi_t)dt.$$ Invoking the finiteness of  $\tau_*$,   we conclude that $\tau_*$ is an optimal stopping time.\smallskip

We now wish to verify that $\tau_*$ coincides with $\tau_D$. Recalling 
  Statement (3) in Lemma \ref{Time}, we proceed to calculate
$$\ba{ll} \ad V\parens{\hat X_{\tau_{D_\e^T}},\hat \varPhi_{\tau_{D_{\e}^{T}}}}\geq V^{T}\parens{0;\hat X_{\tau_{D_{\e}^{T}}},\hat \varPhi_{\tau_{D_{\e}^{T}}}}-\gamma_{T}\parens{\hat \varPhi_{\tau_{D_{\e}^{T}}}} \\
\ns\ad=V^{T}\parens{\tau_{D_\e^T};\hat X_{\tau_{D_{\e_n}^{T}}},\hat \varPhi_{\tau_{D_{\e}^{T}}}}-V^{T}\parens{\tau_{D_\e^T};\hat X_{\tau_{D_{\e}^{T}}},\hat \varPhi_{\tau_{D_{\e}^{T}}}}+V^{T}\parens{0;\hat X_{\tau_{D_{\e}^{T}}},\hat \varPhi_{\tau_{D_{\e}^{T}}}}-\gamma_{T}\parens{\hat \varPhi_{\tau_{D_{\e}^{T}}}} \\
\ns\ad\geq  G\parens{\hat \varPhi_{\tau_{D_{\e}^{T}}}}-\e-\gamma_{T}\parens{\hat \varPhi_{\tau_{D_{\e}^{T}}}}+V^{T}\parens{0;\hat X_{\tau_{D_{\e}^{T}}},\hat \varPhi_{\tau_{D_{\e}^{T}}}}-V^{T}\parens{\tau_{D_\e^T};\hat X_{\tau_{D_{\e}^{T}}},\hat \varPhi_{\tau_{D_{\e}^{T}}}}.\ea  $$
Note that $\tau_*<\infty$ almost surely and that $V^T(t,x,\varphi)$ is bounded and converges to $V(x,\varphi)$ for any fixed $t$ as $T\rightarrow\infty$.
Letting both $T\rightarrow\infty$ and $\e\rightarrow 0^+$, we have $ V(\hat X_{\tau_*},\hat \varPhi_{\tau_*})\geq G(\hat \varPhi_{\tau_*})$. This means that $ \tau_D\leq \tau_*$ almost surely and thus $\tau_*=\tau_D$ almost surely.\smallskip

(4) {\it Optimal stopping boundary}. We now turn to a study of the properties of the optimal stopping boundary. We begin with Lemma \ref{lemmaVN} below.
\begin{lemma}\label{lemmaVN} The following two statements hold.


	{\rm (1)}. For each fixed $x\in\cX$, $ V(x,\cdot)$ is   concave on $[0,\infty)$.\smallskip

	{\rm (2)}. If $x\leq y $ and $\varphi\leq\psi$, $V(x,\varphi)\leq V(y,\psi)$.
	
\end{lemma}
\begin{proof}

	(1) For each fixed $x\in\cX$, $\hat \varPhi_t$ is linear  with respect to the initial value $\varphi$. Further, $ H(x,\cdot)$ and $G(\cdot)$ are concave with respect to $\varphi$ and  $V(x,\cdot)$ is  concave on $(0,\infty)$. The desired result then follows. \smallskip

	(2)  
	Note that $\rho^{-2}(\cdot)>0$ is  increasing on $\cX$. Since $ H(x,\varphi)$ is increasing with respect to $x$ and $\varphi$, the result holds by invoking the comparison of the solutions of stochastic differential equations in \cite[Theorem 1]{Fe2000}. 
\end{proof}

We continue with Proposition \ref{1inC} below.

\begin{proposition}\label{1inC}  The line $L=\cX\times\{1\}$ belongs to the continuation region $  C$. It then follows that $l_0^*(\cdot)<1<l_1^*(\cdot)$.
\end{proposition}
\begin{proof}
	
	Recall that $G(\varphi)=c(1 \wedge\varphi)$. Employing the change-of-variable formula for $t\mapsto G(\hat\varPhi_t)$ and invoking the optional sampling theorem, we have
	$$\EE_{\varphi}G(\hat\varPhi_\t)=G( \varphi)-\frac 12\ell^L_\tau,$$
	where $\ell^L_t$ is the local time of $(\hat X,\hat \varPhi)$ on the curve $c$, i.e.
	$$\ba{ll}\ad\ell^L_t=\PP-\lim_{\e\rightarrow0^+}\frac1{2\e}\int_0^t I\bracks{\hat \varPhi_t\in[1-\e,1+\e]}d\langle\hat\varPhi,\hat\varPhi\rangle_t.\ea$$
	We now define $$ \t_j^t:=t\wedge \inf\{s\geq 0:(\hat X_s,\hat\varPhi_s)\notin (\underline x_j,\overline x_j)\times(1/j,j)\},$$ where $\underline x_j\downarrow\inf\cX$ and $\overline x_j\uparrow\sup\cX$.
	We claim that there exists $t_0>0$ and $\kappa>0$ (independent of $t_0$) such that, for  $t\in(0,t_0)$,
	\bel{localtime}\EE_{x,\varphi}\ell_{\t_K^t}^{L}>\kappa \EE_{x,\varphi}\sqrt{\t_j^t}.\eel
	Assuming the claim in \eqref{localtime} holds, we would then have that for any $x\in\cX$,
	$$ V(x,1)\leq (2K+1)\EE_{x,1} \t^t_j-\k \EE_{x,\varphi}\sqrt{ \t^t_j}+1.$$
	For small values of $t$, since $\PP_{x,\varphi}(\t^t_j>0)=1$, $ V(x,1)<G(1)$. This means that the line $L$ belongs to the continuation region $C$. Now all that remains is to prove \eqref{localtime} as claimed. However, the proof of \eqref{localtime} is exactly the same as that of (8.6) in \cite{ErPZ2020}, and is thus omitted. This completes the proof.
\end{proof}

The stopping region $D$  is divided into two separate parts:
$D_0=D\cap \big(\cX\times (0,1)\big) $ and $ D_1=D\cap \big(\cX\times (1,\infty)\big) $.
Recall the definition of $(l_0^*(\cdot),l_1^*(\cdot))$ in \eqref{l0l1}. The monotonicity of the optimal stopping boundary is now proven in Proposition \ref{PD1} below.

\begin{proposition}\label{PD1}
 The optimal stopping boundary	 $l^*_0(\cdot)$ is increasing in $\cX$. The optimal stopping boundary	 $l^*_1(\cdot)$ is decreasing in $\cX$.


\end{proposition}
\begin{proof} 
Note that for $x_1<x_2$, we have
$$1=V(x_1,l_1^*(x_1))=V(x_2,l_1^*(x_2))\geq V(x_1,l_1^*(x_2)).$$
By the monotonicity of $V(x,\cdot)$, we have that $l_1^*(x_1)\geq l_1^*(x_2)$. This means that $l_1^*(\cdot)$ is decreasing.\smallskip

If $\hat \varPhi_0=0$, then $\hat \varPhi_t=0$ for all $t\geq 0$. The optimal stopping problem of interest then becomes
	$$  V(x,0)=\inf_{\t\in\cT_{\hat X,\hat \varPhi}}\EE_{x,0}\bracks{\int_0^\tau e^{-\lambda t} H( \hat X_t,0)dt}.$$
	Since $ H(x,0)>0$ on $x\in \cX$, instantaneous  stopping 
	is optimal. Note that $V(x,\cdot)$ on $[0,\infty)$ is concave and, for $x\in\cX$, $  V(x,0)=G(0)$.  Then for any $(x,\varphi)$ such that $\varphi< l_0(x)$, we have that $   V(x,\varphi)\geq G(\varphi)$, i.e. $(x,\varphi)\in  D_0$.  Together with the definition of $D$, we have that if $(x,\varphi)\in  D_0$,  then the pairs $(x_0,\varphi_0)$ (with $x_0\geq x$ and $\varphi\geq\varphi_0$) belong to $ D_0$ as well. Therefore, $l^*_0(\cdot)$ must be increasing. This concludes the proof.
\end{proof}

\indent We have now finished the proof of Statement (I) of Theorem \ref{ospt}. We now continue with the proof of  Statement (II) of Theorem \ref{ospt}.\\

\noindent (3) {\it Probabilistic regularity of the optimal stopping boundary}. We need only show Statement (II) of Theorem \ref{ospt} for the case $\varphi_0=1$, as the proof for all other cases is completely identical.	Recall that
	\bel{SDENe}\left\{\ba{ll}\ad d\hat X_t=\hat \mu(\hat X_t)dt+\hat \sigma(\hat X_t)d\hat B_t\\
	\ns\ad  d\hat\varPhi_t=\hat\varPhi_t d\hat B_t.\ea\right.\eel
	For some $x_0\in\cX$, let
	\bel{FU}F(x)=\int_{x_0}^x\frac {\rho(y)}{\sigma(y)}dy,\eel
 and $$\hat U_t=F(\hat X_t)-\log \hat\varPhi_t.$$
	Note that $F$ is continuous and strictly increasing. It\^o's formula then implies that
	\bel{ItoF}\ba{ll}dF(\hat X_t)\ad=F'(\hat X_t)dX_t+\frac12 F''(\hat X_t)\hat\sigma^2(\hat X_t)dt+d\hat B_t\\
	\ns\ad=\parens{\frac{\hat\mu(\hat X_t)}{\hat\sigma(\hat X_t)}-\frac12\hat\sigma'(\hat X_t)}dt+d\hat B_t=K(\hat X_t)dt+d\hat B_t.\ea\eel
By the monotonicity of $b_1$, for any $(x,\varphi)\in S_1$, it is straightforward to see that $(x',\varphi')$ belongs to $D^\circ_1$ if $x'>x$ and $\varphi'>\varphi$.   It then follows from \eqref{ItoF} that, for any $(x,\varphi)\in S_1$,
	\begin{equation}\label{prregular1}\ba{ll}\ad\PP_{x,\varphi}(\t_{ D_1^\circ}<t)\\
	\ns\ad=\PP_{x,\varphi}\((\hat X_s,\hat\varPhi_s)\in  D_1^\circ, \text{ for some }s\in (0,t)\Big)\\
	\ns\ad\geq \PP\(\hat X_s>x,\hat\varPhi_s>\varphi, \text{ for some }s\in (0,t)\Big)\\
	\ns\ad=\PP_{x,\varphi}\(F(\hat X_s)>F(x),\,\log\hat\varPhi_s>\log\varphi,\, \text{for some }s\in (0,t)\Big)\\
	\ns\ad=\PP_{x,\varphi}\parens{\hat B_s+\int_0^sK(\hat X_r)dr>0,\,\hat B_s-\frac s2>0, \text{ for some }s\in (0,t)}
	\\
	\ns\ad=\PP_{x,\varphi}\parens{\frac{\hat B_s}s+\frac 1s\int_0^sK(\hat X_r)dr>0,\frac{\hat B_s}s-\frac 12>0, \text{ for some }s\in (0,t)}=1.\ea\end{equation}
	For any $(x,\varphi)\in S_0$,  we have that for small $\e\in\parens{0,K(x)+\frac12}$, there exists $s_n\downarrow 0$ such that $B_{s_n}=\frac{(1-\e)}2 s_n$. We may thus write
	$$\lim_{n\rightarrow\infty}\frac1{s_n}\int_0^{s_n}K(\hat X_r)dr=K(x)>-\frac12+\e.$$
	Noting that  $(x',\varphi')\in D^\circ_0$ if $x'>x$ and $\varphi'<\varphi$ with $(x,\varphi)\in S_0$, it follows that
	\begin{equation}\label{prregular2}\ba{ll}\ad\PP_{x,\varphi}(\t_{ D_0^\circ}<t)\\
	\ns\ad=\PP_{x,\varphi}\((\hat X_s,\hat\varPhi_s)\in D_0^\circ, \text{ for some }s\in (0,t)\Big)\\
	\ns\ad\geq \PP\(\hat X_s>x,\hat\varPhi_s<\varphi, \text{ for some }s\in (0,t)\Big)\\
	\ns\ad= \PP_{x,\varphi}\(F(\hat X_s)>F(x),\log\hat\varPhi_s<\log\varphi, \text{ for some }s\in (0,t)\Big)\\
	\ns\ad=\PP_{x,\varphi}\(\hat B_s+\int_0^sK(\hat X_r)dr>0,\hat B_s-\frac s2<0, \text{ for some }s\in (0,t)\Big)=1.\ea\end{equation}
	We now combine \eqref{prregular1}  and \eqref{prregular2}. Since $t>0$ was assumed to be arbitrary, it follows that, for any $(x,\varphi)\in S$,
	$\PP_{x,\varphi}(\t_{ D^\circ}=0)=1$. This yields probabilistic regularity of the optimal stopping boundary, as desired.\\

\indent We finish the Appendix with the proof of Statement (III) of Theorem \ref{ospt}.
 It is straightforward to note that  $$\EE_{x_0}\int_0^{\tau^*(\varphi_0)}(1+L_t)dt\leq C\EE_{x_0}\int_0^{\tau^*(\varphi_0)}(1+\varphi_0L_t)dt\leq  C\EE_{x_0,\varphi_0}\tau^*<\infty,$$
	where $\tau^*$ is the optimal stopping time for the optimal stopping problem in \eqref{osp} with initial value $(x_0,\varphi_0)$.
This completes the proof of Theorem \ref{ospt}. \quad \quad \quad \quad \quad \quad \quad \quad \quad \quad \quad \quad \quad \quad \quad \quad \quad \quad \quad \quad $\Box$

\end{document}